\magnification=\magstep1

\def\d{\delta}

\def\im{\item}

\def\ov{\overline}
\def\Ph{\Phi_{\delta}}
\def\moreproclaim{\par}

\def\modo#1{\left|#1\right|}
\def\normo#1{\left\|#1\right\|}

\centerline{\bf Tangent Sequences in Orlicz and Rearrangement Invariant
Spaces}

\bigskip

\centerline{BY PAWE\L\  HITCZENKO}

\medskip
\centerline{\it Department of Mathematics, Box 8205, North Carolina 
State University,}
\centerline{\it Raleigh, NC 27695 -- 8205, USA}

\bigskip
\centerline {AND STEPHEN J.~MONTGOMERY-SMITH}

\medskip
\centerline{\it Department of Mathematics,  University of Columbia -- 
Missouri,}
\centerline{\it Columbia, MO 65211, USA}

\beginsection Abstract

Let $(f_n)$\ and $(g_n)$\ be two sequences of random variables adapted 
to an increasing sequence of $\sigma$-algebras $({\cal F}_n)$\ such that
the conditional distributions of $f_n$\ and $g_n$\ given ${\cal
F}_{n-1}$\ coincide.  Suppose further that the sequence $(g_n)$\ is
conditionally independent. Then it is known that  $\normo{\sum f_k}_p
\le C \, \normo{\sum g_k}_p$, $1 \le p \le
\infty$, where the number $C$\ is a universal constant.
The aim of this paper is to extend this result to
certain classes of Orlicz and rearrangement invariant spaces. This paper
includes fairly general techniques for obtaining rearrangement invariant
inequalities from Orlicz norm inequalities.

\beginsection 1. Introduction

Let $({\cal F}_n)$ be an increasing sequence of $\sigma$-algebras on some
probability space $(\Omega ,{\cal F},P)$. We will assume that ${\cal
F}_0=\{\emptyset ,\Omega\}$. A sequence $(f_n)$ of random variables is
called {\it $({\cal F}_n)$-adapted\/}
if $f_n$ is ${\cal F}_n$-measurable for each $n \ge 1$. In the sequel we
will simply write `adapted' if there is no risk of confusion. For any 
sequence $(f_n)$ of random variables, we will write $f^*=\sup_n|f_n|$ and
$f_n^*=\max_{1\leq k\leq n}|f_k|$. Throughout the paper all equalities
or inequalities between random variables are assumed to hold almost
surely.

Given a $\sigma$-algebra ${\cal A}\subset{\cal F}$ and an integrable random
variable $f$, we
will denote the conditional expectation of $f$ given ${\cal A}$\ by
$E_{\cal
A}f$. If ${\cal A}={\cal F}_k$ then we will simply write $E_kf$ for
$E_{{\cal
F}_k}f$.

The {\it conditional
distribution\/} of a random variable $f$ given ${\cal A}$\ is denoted 
by ${\cal L}(f\big |{\cal A})$. Thus, ${\cal L}(f\big |{\cal A})= {\cal
L}(g\big |{\cal A})$ means that for each real number $t$ we have that
$P(f>t\big | {\cal A})=P(g>t\big | {\cal A})$. 

The following definition was
introduced by Kwapie\'n and Woyczy\'nski in a preprint of their paper
[{\bf 10}], which was distributed as early as 1986. We refer the reader to
their book [{\bf 11}] for more
information on tangent sequences.

\proclaim Definition. {\rm Let $({\cal F}_n)$ be an increasing sequence of
$\sigma$-algebras on $(\Omega ,{\cal F},P)$. \im{(a)} Two adapted sequences
$(f_n)$ and $(g_n)$ of random variables are {\it tangent\/} if for each
$n\ge 1$ we have
$${\cal L}(f_n\big |{\cal F}_{n-1})=
{\cal L}(g_n\big |{\cal F}_{n-1}).$$
\im{(b)} An adapted sequence $(g_n)$ of random variables satisfies {\it
condition~(CI)\/} if there exists a $\sigma$-algebra ${\cal G}\subset {\cal
F}$ such that for each $n\ge1$ $${\cal L}(g_n\big |{\cal F}_{n-1})= {\cal
L}(g_n\big |{\cal G}),$$
and $(g_n)$ is a sequence of ${\cal G}$-conditionally independent
random variables.\moreproclaim}

\noindent A sequence $(f_n)$ is {\it conditionally symmetric\/} if $(f_n)$
and $(-f_n)$ are tangent sequences of random variables. 

Every sequence of random variables $(f_n)$ admits (possibly on an enlarged
probability space) a tangent sequence which satisfies condition~(CI) (cf.
e.g. [{\bf 11} p. 104]). Throughout this paper,
such a sequence will be denoted by $(\overline{f}_n)$, and will be called a
{\it decoupled version\/} of $(f_n)$. It is useful to note that the
$\sigma$-algebra ${\cal G}$ can be chosen so that the random variables
$f_n$, $n \ge 1$, are ${\cal G}$-measurable. 

In this paper, we will be interested in comparing Orlicz and 
rearrangement invariant norms of sums of tangent sequences. A {\it
rearrangement invariant space} $X$ is a space of random variables $f$\
equipped with a complete quasi-norm $\normo {\,\cdot\,}_X$ such that 
either the conditions (i), (ii) and (iii) or (i), (ii) and (iii$'$) below 
 hold: \item{(i)} if $g^{\#} \le f^{\#}$\ and $f \in X$, then $g
\in X$\ with $\normo g_X \le \normo f_X$; \item{(ii)} if $f$\ is simple
with finite support then $f \in X$; \item{(iii)}  $f_n \in X$\ and
$f_n \searrow 0$\ implies $\normo{f_n}_X \searrow 0$ \item{(iii$'$)} 
$f_n \in X$\ and $0 \le f_n \nearrow f$\ and $\sup_n \normo{f_n}_X <
\infty$\ imply $f \in X$\ with $\normo f_X = \sup_n \normo{f_n}_X$.

\noindent
Here $f^{\#}$\ denotes the decreasing rearrangement of $\modo f$, that is,
$f^{\#}(s) = \sup\{ \,t : P( \modo f > t) > s \, \}$. 

Examples of rearrangement invariant spaces include Orlicz spaces and
Lorentz spaces. Let  $\Phi:[0,\infty) \to [0,\infty)$\ be an 
increasing 
function such that $\Phi(0) = 0$, and such that there is
a constant $c>0$\ such that $\Phi(ct) \ge 2 \Phi(t)$\ for all 
$t \ge 0$. (Functions satisfying the latter condition have been 
called {\it dilatory} in [{\bf 14}]; let us note that if $\Phi$ is 
convex, and $\Phi(0) = 0$  then 
this condition is satisfied with $c=2$.) Given such a function 
$\Phi$ we
define the {\it Orlicz norm\/} of a random variable $f$ to be $$ \normo
f_\Phi =
\inf\left\{ \lambda>0 : E(\Phi(\modo{f}/ \lambda)) \le 1 \right\} .$$ We
let
$L_\Phi = \{ f : \normo f_\Phi < \infty\}$. Note that if $\Phi$\ is convex,
then
$L_\Phi$\ is a normed space. However, we do not wish to restrict ourselves
to
normed spaces.

Other examples are the Lorentz spaces. Given $0<p,q \le \infty$,
we define the space $L_{p,q}$\ to be those
random
variables $f$ for which the following quantity is finite: 
$$ \normo f_{p,q}
= \cases{ \displaystyle{
\left({q\over p}\int_0^1 s^{(q/p)-1} f^{\#}(s)^q \, ds \right)^{1/q} } & if $q < \infty$\cr
&\cr
\displaystyle{\sup_{s>0} s^{1/p} f^\#(s) } & if $q = \infty$ .\cr} $$ 
Note that
$L_{p,q}$\ is not a normed space unless $1 \le q \le p \le \infty$. 

Note that the $L_p$\ spaces are special cases: $L_p = L_\Phi =
L_{p,p}$, where $\Phi(t) = t^p$. We refer the reader to [{\bf 13}] for more 
details about these spaces.

In the present paper we will be interested in the domination of a 
rearrangement invariant norm of a sum of an arbitrary sequence of 
adapted random variables  by the rearrangement
invariant norm of a sum of its decoupled version. It is already known
 (see [{\bf 4}]) that if $\Phi$\ satisfies the $\Delta_2$-condition, that is, there is a
constant $c>0$\ such that $\Phi(2t) \le c\, \Phi(t)$\ for all $t \ge 0$,
then there is a constant $C_\Phi$\ such that for every adapted sequence
$(f_n)$ of random variables one has: $$ \|\sum f_i\|_\Phi \leq C_\Phi
\|\sum\overline f_i\|_\Phi . \eqno(1.1)$$ Building on some special
situations considered by Klass [{\bf 4}, Theorem 3.1] and Kwapie\'n 
 [{\bf 9}],
Hitczenko [{\bf 6}] began to investigate how the constant $C_\Phi$\ depends
upon $\Phi$. He showed that there is a universal constant $C>0$\ such
that $$ \|\sum f_i\|_p \leq C\, \|\sum\overline f_i\|_p \qquad 1 \le p
\le \infty . \eqno(1.2)$$ In this present paper, we will  show among
other things, that inequality $(1.1)$\ holds with $C_\Phi$\ uniformly
bounded, at least for certain classes of Orlicz functions. 
Our first theorem extends a result of Klass who proved (1.1) for
randomly stopped sums of independent random variables.  

Let us define classes of Orlicz functions. Following Klass [{\bf 8}], for
$q>0$ we define the class $F_q$ as the class of all functions $\Phi$ such
that: \item{(i)}
$\Phi:[0,\infty)\rightarrow[0,\infty)$, $\Phi(0)=0$, \item{(ii)} $\Phi$ is
nondecreasing and continuous, and \item{(iii)} $\Phi$ satisfies the growth
condition: $\Phi(cx)\le c^q\Phi(x)$ for all $x\ge 0$, $c\ge2$. 

\noindent For $p >0$,
we define the class $G_p$\ as the class of all functions $\Phi$ such that
\item{(i)} $\Phi:[0,\infty)\rightarrow[0,\infty)$, $\Phi(0)=0$, \item{(ii)}
$\Phi$ is nondecreasing and continuous, and \item{(iii)} $\Phi$ satisfies
the growth condition: $\Phi(cx)\ge c^p\Phi(x)$ for all $x\ge 0$, $c\ge2$. 

Then we obtain the following results.

\proclaim Theorem 1.1. There is a universal constant $C>0$\ such that if
$\Phi \in F_q$\ for some $q > 0$, then for every adapted sequence
$(f_n)$ of random variables one has:  $$ E \Phi\left(\modo{\sum
f_i}\right) \leq C^{1+q} E \Phi\left(\modo{\sum \overline f_i}\right).$$ 

This inequality had already been obtained by Klass [{\bf 8}] in the special
case
that $f_k=I(\tau\ge k)\xi_k$, where $(\xi_k)$ is a sequence  of
independent random variables and $\tau$ is a stopping time. More
precisely,  Klass proved his result for Banach space valued random
variables $(\xi_k)$ (with absolute value replaced by norm). To discuss 
Banach space valued random variables one needs to adjust notation; 
for a random
variable $Y$ and a $\sigma$-algebra ${\cal A}$, we use 
${\cal L}(Y |{\cal A})$ to denote the regular version of the conditional 
distribution of $Y$ given ${\cal A}$, that is, ${\cal L}(Y |{\cal A}) = 
{\cal L}(Z |{\cal A})$ means that for every Borel subset $A$\ of the 
Banach space, we 
have
$P(Y\in A|{\cal A}) = P(Z\in A|{\cal A})$. Recall that the existence of the 
regular versions of the conditional distributions is guaranteed, 
as long as our random variables take values in a separable Banach space.
As it turns out,
in our generality, the inequality of Theorem~1.1 need not hold (with
{\it any\/} constant), unless some extra  conditions are imposed 
on the geometry of the
underlying Banach space (see e.g. [{\bf 3}]). Since it is unclear
at this time for which Banach spaces the inequality $$
\big(E\|\sum f_k\|^p\big)^{1/p}\le c_p\big(E\|\sum
\ov{f}_k\|^p\big)^{1/p} $$ holds (even if the constant $c_p$ is allowed
to depend on $p$), we confine our discussion to real valued random
variables. 

\proclaim Corollary 1.2. Given numbers $p_0>0$\ and $r\ge 1$,
there is a constant $c_{p_0,r}$\ such that if $p \ge p_0$, and if $\Phi \in
G_p \cap F_{rp}$, then $$ \|\sum f_i\|_\Phi \leq c_{p_0,r} \|\sum\overline
f_i\|_\Phi.$$

The next step is to extend these results to rearrangement invariant 
spaces.  This will be accomplished through a rather general
method of obtaining rearrangement invariant norm inequalities from Orlicz
norm inequalities. We believe that this technique will prove useful
in other contexts as well. In particular, we would like to mention
that this method could be used to deduce martingale inequalities
obtained by Johnson and Schechtman [{\bf 7}] from the corresponding
inequalities for Orlicz functions.

   Corresponding to the notions of
Orlicz spaces lying in $G_p \cap F_q$, we have the following notion. We
say that a rearrangement invariant space is an {\it interpolation space
for $(L_p,L_q)$\/} (in short, a {\it $(p,q)$-interpolation space\/}) if
there is a constant $c > 0$ such that for every operator $T:L_p \cap L_q
\to L_p \cap L_q$\ for which $\normo T_{L_p \to L_p} \le 1$\ and $\normo
T_{L_q \to L_q} \le 1$\ we have that $\normo T_{X \to X} \le c$.

However, this notion is not quite what we need. Define $$ K_{p,q}(f,t) =
\inf\{ \normo{f'}_p + t \normo{f''}_q : f' + f'' = f^\# \} .$$
We will say that a rearrangement invariant space $X$\ is a  {\it
$(p,q)$-$K$-interpolation space\/} if there is a constant 
$c$\ such that
whenever $f$\ and $g$\ are such that 
$K_{p,q}(f,t) \le K_{p,q}(g,t)$\ ($t > 0$),   and $g \in
X$, then $f \in X$\ and $\normo f_X \le c \normo g_X$. The {\it
$(p,q)$-$K$-interpolation constant\/} of $X$, denoted by $C_{p,q}(X)$,
is the infimum of $c$\ that work for all  functions $f$\ and $g$.

It is quite easily seen that every $(p,q)$-$K$-interpolation space is a
$(p,q)$-interpolation space.
It is also known that if $1 \le p,q \le \infty$, then every normed
$(p,q)$-interpolation space is a 
$(p,q)$-$K$-interpolation space (see [{\bf 1}]).
We are able to establish the following method for obtaining 
rearrangement invariant inequalities from Orlicz inequalities.

\proclaim Theorem 1.3. Suppose that $\normo f_\Phi \le  \normo g_\Phi$\
for all $\Phi \in F_q \cap G_p$, where $0 < p < q < \infty$. If
$X$\ is a $(p,q)$-$K$-interpolation space, then $\normo f_X \le
2^{2+1/p} C_{p,q}(X) \normo g_X$.

\proclaim Corollary 1.4. Given numbers $p_0>0$\ and $r\ge 1$,
there is a constant $c_{p_0,r}$\ such that if $p \ge p_0$, and if $X$\ is a
$(p,pr)$-$K$-interpolation space, then $$ \|\sum f_i\|_X \leq c_{p_0,r}
C_{p,pr}(X) \|\sum\overline f_i\|_X.$$

In particular, we are able to obtain the following result for  Lorentz
spaces.  Note that if one is interested in normed Lorentz
spaces, then $p_0$ below can be taken to be 1, and
the resulting inequality extends (1.2).

\proclaim Corollary 1.5. Given a number $p_0>0$, there is a constant
$c_{p_0}$\ such that if $p,q \ge p_0$, then $$ \|\sum f_i\|_{p,q} \leq
c_{p_0}
\|\sum\overline f_i\|_{p,q}.$$

\beginsection 2. Inequalities for Orlicz functions 

We begin with a proof of Theorem~1.1. Since our proof is based on 
well understood techniques we will be somewhat sketchy and we refer the
reader to [{\bf 6}] for details that are not explained here.
Throughout this section we let $(M_n)$ be a martingale with difference
sequence $(\Delta_k)$. Since $F_{q_1}\subset F_{q_2}$ whenever $q_1\le
q_2$, we can assume without loss of generality that $q\ge 1$. Our
departing point is the following result which can be found in the just
mentioned paper (Theorem~5.1 and the beginning of the proof of
Lemma~2.3). 

\proclaim Lemma 2.1. Let $1\le q<\infty$, and let $(\Delta_k)$ be a
conditionally symmetric martingale difference sequence, and $(\ov{\Delta
}_k)$ its decoupled version. Set $T_{n,q}(M)=(E|\sum_{k=1}^n\ov
{\Delta}_k|^q\big|{\cal G})^{1/q}$. Then there exist $\delta_1>0$,
$\beta>1+\delta_1$ and $\epsilon$\ with
$0<\epsilon\le 1/2$ such that for every
$\lambda>0$ we have
$$ P(M^*\geq \beta\lambda, (T_q^*(M)\vee \Delta^*)<\delta_1\lambda) \leq
\epsilon^qP(M^*\ge\lambda). $$

\noindent From this, we obtain

\proclaim Lemma 2.2. Let
$(\Delta_k)$ be as above, and assume that $w_n$ is a ${\cal
F}_{n-1}$-measurable random variable such that $|\Delta_n|\le w_n$\ 
for each $n \ge 1$.
Set $N_n=\sum_{k=1}^n\ov {\Delta}_k$ ($n\ge 1$). Suppose that
$\delta_1$, $\beta$, and $\epsilon$ are as in Lemma~2.1. Then, there
exist $\delta>0$, $\delta_2>0$ and $0<\alpha\le 1/2$ such that for every
$\lambda>0$ we have $$ P(M^*\geq \beta\lambda, N^*<\delta\lambda) \leq
\epsilon^qP(M^*\ge\lambda) +P(w^*\ge\delta_2\lambda) +
(1-\alpha^q)P(M^*\geq \beta\lambda). $$

\noindent{\bf Proof:}
We have that
$$
\eqalign{P(M^*\ge\beta\lambda,N^*<
\delta\lambda)\le& P(M^*\ge\beta\lambda,T_q^*(M)< \delta_1\lambda, w^*<
\delta_2\lambda)+P(w^*\ge\delta_2\lambda)\cr+&
P(M^*\ge\beta\lambda,T_q^*(M)\ge
\delta_1\lambda,w^*<\delta_2\lambda,N^*< \delta\lambda). \cr}\eqno(2.1) $$
Suppose $\delta_2\le\delta_1$. Then, in view of Lemma 2.1, for the first
probability on the right-hand side of (2.1) we have that $$
\eqalign{P(M^*\ge\beta\lambda,T_q^*(M)< \delta_1\lambda, w^*<
\delta_2\lambda)\le&
P(M^*\ge\beta\lambda,T_q^*(M)\vee w^*< \delta_1\lambda)\cr\le&
\epsilon^qP(M^*\ge\lambda).\cr} $$
It remains to estimate the last probability in (2.1). Since $M^*$, $w^*$
and
$T_q^*(M)$ are ${\cal G}$-measurable, by conditioning on ${\cal G}$, we see
that
the last probability in (2.1) is equal to: $$
E\Big\{I(M^*\ge\beta\lambda,T_q^*(M)\ge
\delta_1\lambda,w^*<\delta_2\lambda)P(N^*< \delta\lambda|{\cal G})\Big\}.
$$ By
 Kolmogorov's converse inequality (see e.g. [{\bf 12},
Remark
6.15, p. 161]) for all sequences of independent and symmetric random 
variables $(\xi_k)$ and for all $t>0$ we have that $$ P(S^*\ge
t)\ge {1\over 2^q}\Big(1-{2^{2q}(t^q+E(\xi^*)^q)\over E(S^*)^q}\Big), $$
where $S_n=\sum_{k=1}^n\xi_k$. Applying this result conditionally on
${\cal G}$, we obtain that $$ P(N^*<\delta\lambda|{\cal G})\le 1-{1\over
2^q}\Big(1-{2^{2q}((\delta\lambda)^q+E_{\cal G}(\ov {\Delta}^*)^q)\over
E_{\cal G}(N^*)^q}\Big).\eqno(2.2) $$ Also, if $w_n<\delta_2\lambda$,
then $|\Delta_n|<\delta_2\lambda$, and since $w_n$ is ${\cal
F}_{n-1}$-measurable, and the conditional distributions of $\Delta_n$
and $\ov {\Delta}_n$ coincide, it follows that $|\ov{\Delta}_n|<
\delta_2\lambda$. Therefore, on the set $$
\{T_q^*(M)\ge\delta_1\lambda,w^*<\delta_2\lambda\}, $$ we have $$
{(\delta\lambda)^q+E_{\cal G}(\ov{ \Delta}^*)^q\over E_{\cal
G}(N^*)^q}\le{(\delta\lambda)^q+(\delta_2\lambda)^q\over
(\delta_1\lambda)^q}, $$ so that the conditional probability in (2.2)
does not exceed $$ 1-{1\over 2^q}\Big(1-{4^q(\delta^q+\delta_2^q)\over
\delta_1^q}\Big). $$ Choosing
$\delta=\delta_2=\delta_1/12$, we obtain that $$ 1-{1\over
2^q}\Big(1-{4^q(\delta^q+\delta_2^q)\over \delta_1^q}\Big)=1-{1\over
2^q}\Big(1-{2\over 3^q}\Big)\le 1-\alpha^q, $$ whenever $\alpha\le 1/6$.
Therefore, $$ \eqalign{ E\Big\{I(M^*\ge\beta\lambda&,T_q(M)\ge
\delta_1\lambda,w^*<\delta_2\lambda)P(N^*< \delta\lambda|{\cal
G})\Big\}\cr\le& (1-\alpha^q) EI(M^*\ge\beta\lambda,T_q(M)\ge
\delta_1\lambda,w^*<\delta_2\lambda)
\le(1-\alpha^q)P(M^*\ge\beta\lambda).\cr} $$ This completes the proof of
Lemma 2.2. 

\bigskip
\noindent Now we are ready to complete the proof of Theorem~1.1. By an argument
similar to one used in [{\bf 6}, proof of Lemma~2.1], it follows
that
in order to prove
$$ E\Phi(|\sum f_k|)\le
c^qE\Phi(|\sum\ov{f}_k|),\eqno(2.3) $$ it suffices to
establish (2.3) for $(f_k)=(\Delta_k)$, a conditionally symmetric
martingale difference sequence. By a routine application of Davis'
decomposition (cf. e.g. [{\bf 2}] and references therein), we
may also assume that $|\Delta_n|\le w_n$, where $w_n$ is a ${\cal
F}_{n-1}$-measurable random variable, and that $w^*\le 2\Delta^*$. The
latter inequality, together with the inequality $P(f^*\ge t)\le
2P(g^*\ge t)$ valid for all tangent sequences $(f_k)$ and $(g_k)$ 
(cf. [{\bf  4}] or 
[{\bf 11}, Theorem 5.2.1 (i)]), implies that $$
P(w^*\ge t)\le P(\Delta^*\ge t/2)\le 2P(\ov{\Delta}^*\ge t/2)\le
2P(N^*\ge t/4).\eqno(2.4) $$ By Lemma 2.2, we have that $$
P(M^*\ge\beta\lambda)\le
P(N^*\ge\delta\lambda)+P(w^*\ge\delta_2\lambda)+
(1-\alpha^q)P(M^*\ge\beta\lambda), $$ so that $$
\alpha^qP(M^*\ge\beta\lambda)\le
P(N^*\ge\delta\lambda)+P(w^*\ge\delta_2\lambda). $$ Consequently, $$
\alpha^qE\Phi(M^*/\beta)\le E\Phi(N^*/\delta)+E\Phi(w^*/\delta_2). $$
Therefore, $$ \eqalign{\Big({\alpha\over\beta}\Big)^qE\Phi(M^*)=&
\Big({\alpha\over\beta}\Big)^qE\Phi(\beta M^*/\beta)\cr
\le&\Big({\alpha\over\beta}\Big)^q\beta^pE\Phi( M^*/\beta)\le
E\Phi(N^*/\delta)+E\Phi(w^*/\delta_2)\cr
\le&\delta^{-q}E\Phi(N^*)+\delta_2^{-q}E\Phi(w^*).\cr} $$ By (2.4) we
have that $$ E\Phi(w^*)\le 2E\Phi(4N^*)\le 2\cdot4^qE\Phi(N^*), $$ so
that $$ E\Phi(M^*)\le
\Big({\beta\over\alpha}\Big)^q\Big\{{1\over\delta^q}
+{2\cdot4^q\over\delta_2^q}\Big\}E\Phi(N^*). $$ 
Conditionally given ${\cal G}$, $(\ov\Delta_k)$ is a sequence of independent and symmetric 
random variables. Therefore, by Levy's
inequality, 
$$P(N^*\ge t\big|{\cal G})\le 2P(N\ge t\big|{\cal G}).$$
This implies (see e.g. [{\bf 11}, Proposition 0.2.1]) that for every increasing function 
$\phi:{\bf R}^+\to{\bf R}^+$\ we have $E\phi(N^*)\le2E\phi(N)$. This
 completes the proof of Theorem 1.1.

As we mentioned in the introduction,
inequality (2.3) extends a result of Klass, who considered sequences
$(f_k)$ of special form. On the other hand, it follows from a result of
Kwapie\'n [{\bf 9}] that if $f_k=(\sum_{j=1}^{k-1}a_{j,k}\xi_j)\xi_k$,
where $(\xi_j)$ is a sequence of independent zero mean random variables,
then (2.3) holds in the stronger form: $$ E\Phi(|\sum f_k|)\le
E\Phi(c|\sum\ov f_k|), $$ for some absolute constant c and for {\it
every} convex function $\Phi$. Thus one may wonder whether our
restrictions on $\Phi$ can be relaxed. We wish to close this section
with a negative result showing that (2.3) does not
hold for all convex functions $\Phi$. (However, it is still possible that (2.3)
holds under weaker assumption than ours.) Our example is an easy
adaptation of an example due to Talagrand concerning comparison of tail
behavior for sums of tangent sequences. This example was included in
[{\bf 5}], and we refer the reader to the latter paper for 
details that are not included here. 

\proclaim
Proposition 2.3. For every constant $c>0$, there exists a convex function
$\Phi:\bf R\to\bf R$ and a sequence $(f_k)$ such that $$ E\Phi(|\sum
f_k|)\ge cE\Phi(c|\sum\overline{f}_k|). $$ 

\noindent{\bf Proof:}\ \
We will show that for every
$k\in{\bf N}$ there exists a convex function $\Phi$ and a sequence $(f_k)$
for which $$ E\Phi(|\sum f_k|)\ge {2^{2^{k+2}}\over
k^22^{2k}}E\Phi({k\over4}|\sum\overline{f}_k|). $$ Let $(r_n)$\ denote the
Rademacher random variables, that is, a sequence of independent random
variables such that $P(r_n = \pm 1) = 1/2$. Fix $k\in{\bf N}$. Given an
integer $N_1$ to be specified in a moment, define $N_2,\dots N_k$ as
follows: $$
N_i-N_{i-1}=2^{-(i-1)}N_1,\quad i=2,\dots ,k.  $$ Put $$ \Omega _1=\{
r_1=\dots =r_{N_1}=1\} , $$ and then $$ \Omega _i=\Omega _{i-1}\cap\{
r_{N_{i-1}+1}=\dots =r_{N_i}\},\quad i=2,\dots ,k. $$ Define a sequence
of random variables $(v_i)$ by the formulas: $$ \eqalign {v_1=\dots
=v_{N_1}&=1\cr v_{N_1+1}=\dots =v_{N_2}&=2I_{\Omega _1}\cr
\dots\dots\dots\dots &\cr v_{N_{k-1}+1}=\dots =v_{N_k}&=2^{k-1}I_{\Omega
_{k-1}}.\cr} $$ We let $f_j=v_jr_j$ for $j=1,\dots N_k$. Then
$\overline{f}_j=v_jr_j'$, where $(r_j')$ is an independent copy of
$(r_j)$ (cf. [{\bf 10}, Example 4.3.1]). For
$0<\delta<1$, let $\Phi_{\delta}$ be a convex function defined by
$\Ph(x)=(x-\delta kN_1)^+$. Note that $\sum |v_j|=kN_1$, and therefore
$$ E\Ph(|\sum v_jr_j|)\ge (1-\d)kN_1P(|\sum v_jr_j|\ge kN_1). $$ On the
other hand, if $|\sum v_jr_j'|<4N_1$, then $|\sum v_jr_j'|\le
4N_1-1$, so that with $\d=1-1/(4N_1)$ we get $$ \Ph({k\over4}|\sum
v_jr_j'|)\le(kN_1-{k\over 4}-\d kN_1)^+=0. $$ Since $$ P(|\sum
v_jr_j'|\ge 4N_1)\le k2^{-N_1/2^{k-2}}P(|\sum v_jr_j|\ge kN_1), $$ (cf.
[{\bf 5}, top half of page 176]) we obtain $$
\eqalign{E\Ph({k\over4}|\sum v_jr_j'|)\le& \Ph({k^2N_1\over4})P(|\sum
v_jr_j'|\ge 4N_1)\cr \le& kN_1{k\over4}k2^{-N_1/2^{k-2}}P(|\sum
v_jr_j|\ge kN_1),\cr} $$ and it follows that $$ {E\Ph(|\sum
v_jr_j|)\over E\Ph((k/4)|\sum v_jr_j'|)}\ge
{4(1-\d)kN_12^{N_1/2^{k-2}}\over k^3N_1}= {2^{N_1/2^{k-2}} \over
k^2N_1}\ge{2^{2^{k+2}}\over k^22^{2k}}, $$ for $N_1\ge 2^{k}$. This
completes the proof.

In the above example the sequence $(f_k)$ may be constructed so that $\sum
f_k$ is a randomly stopped sum of independent random variables (see Remark
on p. 176 of [{\bf 5}]). Thus, the conclusion of this Remark applies
here as well.

\beginsection 3. Rearrangement invariant norm inequalities 

\proclaim Lemma 3.1. Suppose that $\Phi$\ and $\Psi$\ are two Orlicz
functions. Let $\Theta = \Phi \wedge \Psi$, and $\Theta_1(x) = {1\over 2}
\Theta(x)$. Then $$
\textstyle{1\over 2}\normo f_{\Theta_1}
\le \inf\{ \normo{f'}_\Phi + \normo{f''}_\Psi : f' + f'' = f^\#\} \le
2\,\normo f_{\Theta} .$$

\noindent{\bf Proof:} To show the left hand side, suppose that  $f^\# =
f' + f''$, and $\normo{f'}_\Phi + \normo{f''}_\Psi \le 1$. Then 
$E\Phi(\modo{f'}) \le 1$\ and $E\Psi(\modo{f''}) \le 1$, and so
$$
\eqalign{ E\Theta_1(\textstyle{1\over 2}|f|) \le& \textstyle{1\over 2}
E\Theta(\max\{\modo{f'},\modo{f''}\}) = \textstyle{1\over 2} E
\max\{\Theta(\modo{f'}),\Theta(\modo{f''})\} \cr
\le& \textstyle {1\over 2}
(E\Phi(|f'|) + E\Psi(|f''|)) \le 1 .\cr}$$ To show the right hand side,
suppose that $\normo f_{\Theta} \le 1$, that is, $E\Theta(|f|) \le 1$.
Let $$ f'(t) = \cases{ f^\#(t) & if $\Phi(|f(t)|) \le \Psi(|f(t)|)$ \cr
0 & otherwise,\cr} $$ and $f'' = f^\# - f'$. Then we see that
$E\Phi(|f'|) \le E\Theta(\modo{f}) \le 1$\ and that $E\Psi(|f''|) \le
E\Theta(\modo{f}) \le 1$, and the result follows. 

The next lemma follows immediately. 

\proclaim
Lemma 3.2. Let $\Phi_t(x) = x^p \wedge (tx)^q$, where $0 < p < q < \infty$.
Then $$ 2^{-1-1/p} \normo f_{\Phi_t} 
\le K_{p,q}(f,t) \le 2 \normo f_{\Phi_t} .$$

Now we will prove Theorem~1.3, using the above Lemma. 
From the hypothesis of Theorem~3.1, it follows that $\normo
f_{\Phi_t} \le \normo g_{\Phi_t}$. Hence by Lemma 3.2, it follows that 
$K_{p,q}(f,t) \le 2^{2+1/p}K_{p,q}(g,t) $. Now the result follows
by the definition of $(p,q)$-$K$-interpolation space. 

Corollary~1.4 follows easily from Theorems 1.1 and 1.3. To show
Corollary~1.5, we only need the following result. The methods below are
all fairly standard in interpolation theory, and indeed if one is not
concerned about uniform estimates, may be taken directly from the
literature. 

\proclaim Lemma 3.3. Given $p_0>0$, there is a constant $c_{p_0}>0$\ such
that if $p,q \ge p_0$, then $L_{p,q}$\ is a 
$(p/2,2p)$-$K$-interpolation space with constant bounded by $c_{p_0}$. 

\noindent {\bf Proof:}\ \ First let us define some norms. For $p\le q$, let
$$ \normo f_{a(t)} = \inf\{ t^{-2/p} \normo{f'}_{p/2} + t^{-1/2p}
\normo{f''}_{2p} : f'+f'' = f^\# \} ,$$
$$ \normo f_{b(t)} = \left( {1\over t} \int_0^t f^\#(s)^{p/2} \, ds
\right)^{2/p}
+ \left( {1\over t} \int_t^\infty
f^\#(s)^{2p} \, ds \right)^{1/2p} .$$
Clearly $\normo f_{a(t)} \le \normo f_{b(t)}$. Also $\normo f_{a(t)} \ge
\min\{1,2^{1-2/p}\} f^\#(t)$.
This is because if $f^\# = f' + f''$, then $$ \eqalignno{ t^{-2/p}
\normo{f'}_{p/2} + t^{-1/2p} \normo{f''}_{2p} &\ge \left({1\over t}
\int_0^t \modo{f'(s)}^{p/2} \, ds \right)^{2/p} + \left({1\over t} \int_0^t
\modo{f''(s)}^{2p} \, ds \right)^{1/2p} \cr &\ge \left({1\over t} \int_0^t
\modo{f'(s)}^{p/2} \, ds \right)^{2/p} + \left({1\over t} \int_0^t
\modo{f''(s)}^{p/2} \, ds \right)^{2/p} \cr &\ge \min\{1,2^{1-2/p}\}
\left({1\over t}
\int_0^t f^\#(s)^{p/2} \, ds \right)^{2/p} \cr &\ge
\min\{1,2^{1-2/p}\} f^\#(t) .\cr }$$
Next, given a function $f$, let us define the function $Hf(t) = \normo
f_{b(t)}$. Then it follows that $\normo{Hf}_{p,q} \le
32^{1/\min\{p,q\}}\normo f_{p,q}$. To see this, first note
that $\normo{Hf}_{p,q} \le (\normo{H_1f}_{p,q}^q +
\normo{H_2f}_{p,q}^q)^{1/q}$,
where $$ \eqalignno{ H_1 f(t) &= \left({1\over t} \int_0^t f^\#(s)^{p/2} \,
ds,
\right)^{2/p} \cr H_2 f(t) &= \left({1\over t} \int_t^\infty f^\#(s)^{2p}
\, ds
\right)^{1/2p} \cr} $$ We will use two properties of $L_{p,q}$. First, if
$v \le
q$, and if $f_1$, $f_2,\dots,$\ $f_n$\ are functions, then $$
\normo{\left(\sum_{i=1}^n (f_i^\#)^v\right)^{1/v} }_{p,q} \le
\left(\sum_{i=1}^n
\normo {f_i}_{p,q}^v \right)^{1/v} .$$ Second, if we define the operators
$D_af(t) = f^\#(at)$\ for $0 < a < \infty$, 
then $\normo{D_af}_{p,q} =
a^{-1/p}
\normo f_{p,q}$.
 (The first property follows from Minkowski's inequality for $L_
{q/v}$, 
the second is a simple change of variables argument.)
Let $u = \max\{p/q,2\}$, $v = \min\{q,p/2\}$. Then $$ \eqalignno{
\normo{H_1f}_{p,q}
&\le
\normo{\left(\int_0^1 (D_a f^\#)^{p/2} \, da \right)^{2/p}}_{p,q} \cr &\le
\normo{\left(\sum_{n=0}^\infty \int_{2^{-u (n+1)}}^{2^{-u n}} (D_a
f^\#)^{p/2} \, da \right)^{2/p} }_{p,q} \cr &\le
\normo{\left(\sum_{n=0}^\infty 2^{-u n}
(D_{2^{-u (n+1)}} f^\#)^{p/2} \right)^{2/p} }_{p,q} \cr &\le
\normo{\left(\sum_{n=0}^\infty 4^{-n} (D_{2^{-u(n+1)}} f^\#)^v
\right)^{1/v} }_{p,q} \cr &\le
\left( \sum_{n=0}^\infty 4^{-n} \normo{D_{2^{-u(n+1)}} f^\#}_{p,q}^v
\right)^{1/v}
\cr
&\le
\left( \sum_{n=0}^\infty 4^{-n} 2^{u(n+1)v/p})\right)^{1/v} \normo f_{p,q}
\cr
&\le
\left( \sum_{n=0}^\infty 2^{1-n} \right)^{1/v} \normo f_{p,q} \cr &\le
4^{1/v} \normo f_{p,q} .\cr}$$
Now let $u = \max\{2p/q,1\}$, $v = \min\{q,2p\}$. $$ \eqalignno{
\normo{H_2f}_{p,q}
&\le
\normo{\left(\int_1^\infty (D_a f^\#)^{2p} \, da \right)^{1/2p}}_{p,q} \cr
&\le
\normo{\left(\sum_{n=0}^\infty \int_{2^{u n}}^{2^{u(n+1)}} (D_a f^\#)^{2p}
\, da \right)^{1/2p} }_{p,q} \cr &\le \normo{\left(\sum_{n=0}^\infty
2^{u(n+1)} (D_{2^{u n}} f^\#)^{2p} \right)^{1/2p} }_{p,q} \cr &\le
\normo{\left(\sum_{n=0}^\infty 2^{n+1} (D_{2^{u n}} f^\#)^v \right)^{1/v}
}_{p,q} \cr &\le
\left( \sum_{n=0}^\infty 2^{n+1}
\normo{D_{2^{u n}} f^\#}_{p,q}^v \right)^{1/v} \cr &\le \left(
\sum_{n=0}^\infty 2^{n+1} 2^{-u n v/p})\right)^{1/v} \normo f_{p,q} \cr
&\le
\left( \sum_{n=0}^\infty 2^{1-n} \right)^{1/v} \normo f_{p,q} \cr &\le
4^{1/v} \normo f_{p,q} .\cr}$$
Finally, to finish, suppose that  $K_{p/2,2p}(f,t) \le K_{p/2,2p}(g,t)$\
for all $t > 0$. Then it follows that  
$$ 
\eqalign{f^\#(t) \le& \max\{1,2^{2/p-1}\}
\normo f_{a(t)} \le \max\{1,2^{2/p-1}\} \normo g_{a(t)} \le 
\max\{1,2^{2/p-1}\} \normo g_{b(t)} \cr
=& 2^{2/p} Hg(t) .\cr}
$$ Hence $$
\normo f_{p,q} \le 2^{2/p} \normo{Hg}_{p,q} \le 128^{1/\min\{p,q\}}
\normo g_{p,q} .$$ 

\beginsection Acknowledgments

The second named
author would like to express warm gratitude to Victor de la Pe\~na  for
introducing him to this problem. Research of both authors was partially
supported by separate NSF grants.  The second named author was also
supported by the
Research Board of the University of Missouri.

\beginsection References

\medskip
\frenchspacing
\newcount\refnum
\refnum=0
\def\ref{\global\advance\refnum by 1\item{[{\bf\the\refnum}]}} 

\ref J. ARAZY and M.  CWIKIEL.  A new characterization of the
interpolation spaces between $L^p$\ and $L^q$, {\it Math. Scand.} {\bf 
55} (1984),
253--270.

\ref  D. L. BURKHOLDER. Distribution function inequalities for
martingales, {\it Ann. Prob\-ab.} {\bf 1} (1973), 19 - 42. 

\ref  D. J. H. GARLING. Random martingale transform inequalities,
{\it
Probability in Banach Spaces, 6 (Sandbjerg, Denmark, 1986),  101 -- 119,
Progr.
Probab.} {\bf 20}, Birkh\"auser, Boston (1990). 

\ref P. HITCZENKO.  Comparison of moments for tangent sequences
 of random variables. {\it  Probab. Theory Related Fields} {\bf 78} 
(1988), 223 -- 230.

\ref P. HITCZENKO.  Domination inequality for martingale transforms
of a Radema\-cher sequence. {\it Israel J. Math.} {\bf 84} (1993), 161 -- 178. 

\ref P. HITCZENKO.  On a domination of sums of random variables by
sums
of conditionally independent ones. {\it Ann. Probab.} {\bf 22} (1994), 453 -- 468.

\ref  W. B. JOHNSON and  G. SCHECHTMAN.  Martingale inequalities in
rearrangement invariant function spaces, {\it Israel J. Math.} 
{\bf 64} (1988),
 267 - 275.

\ref  M. J. KLASS. 
 A best possible improvement of
Wald's equation, {\it   Ann. Probab.} {\bf 16} (1988), 840 - 853.

\ref S. KWAPIE\'N.
Decoupling inequalities for polynomial chaos, {\it Ann. Probab.} {\bf 15}
 (1987), 
1062 - 1072.

\ref S. KWAPIE\'N and W. A. WOYCZY\'NSKI.  Semimartingale integrals
via decoupling inequalities and tangent processes, {\it Probab. Math.
Statist.} {\bf 12} (1991), 165 -- 200.

\ref  S. KWAPIE\'N and W. A. WOYCZY\'NSKI.  {\it Random Series and
Stochastic Integrals. Single and Multiple}, Birkh\"auser, Boston (1992).

\ref M. LEDOUX and M. TALAGRAND.  {\it Probability in Banach
Spaces.} Springer, Berlin, Heidelberg (1991).

\ref J. LINDENSTRAUSS and L.  TZAFRIRI.  {\it Classical Banach Spaces. Function Spaces}, Springer, Berlin, Heidelberg, New York (1977). 

\ref S. J. MONTGOMERY -- SMITH. Comparison of Orlicz -- Lorentz spaces,
{\it Studia Math.} {\bf 103} (1992), 161 -- 189.

\bye